\documentclass[12pt]{article}

 \newcommand{\Title}{Geometric Ergodicity and Perfect Simulation}
 
 \def\Keywords{CFTP, domCFTP, geometric Foster-Lyapunov condition, 
geometric ergodicity, Markov chain Monte Carlo, perfect simulation, 
uniform ergodicity}

 \usepackage{ifthen,amsmath,amsfonts,multicol,calc,times}
 \usepackage{chicago,url,xspace}

\usepackage{hyperref}

 \newtheorem{result}{Result}
 \newtheorem{algorithm}[result]{Algorithm}
 \newtheorem{definition}[result]{Definition}
 \newtheorem{theorem}[result]{Theorem}
 \newtheorem{lemma}[result]{Lemma}
 \newenvironment{proof}[1][]{\noindent%
\ifthenelse{\equal{#1}{}}{\textbf{Proof:} }%
{\textbf{Proof (#1):\newline}}}%
{\smallskip\hfill\(\Box\)\bigskip}

 \renewcommand{\d}{\operatorname{\text{d}}}
 \newcommand{\dist}{\operatorname{\text{dist}}}
 \newcommand{\CFTP}{\emph{CFTP}\xspace}
 \newcommand{\domCFTP}{\emph{dom}\CFTP}
 \newcommand{\Expect}[1]{\operatorname{\mathbb{E}}\left[#1\right]}
 \newcommand{\Indicator}[1]{\operatorname{\mathbb{I}}\left[#1\right]}
 \newcommand{\Law}[1]{\mathcal{L}\left({#1}\right)}
 \newcommand{\Leb}{\operatorname{\text{Leb}}}
 \newcommand{\Prob}[1]{\operatorname{\mathbb{P}}\left[#1\right]}

\begin{document}


\thispagestyle{empty}

\title{\Title}
\author{Wilfrid S.~Kendall}

\date{\text{ }}

\maketitle

\begin{quote}
\end{quote}
\begin{quote}
\end{quote}

\begin{abstract}
  This note extends the work of \citeN{FossTweedie-1997}, who showed
  that availability of the classic \citeN{ProppWilson-1996} 
Coupling from The Past algorithm
  is essentially equivalent to uniform
  ergodicity for a Markov chain (see also \citeNP{HobertRobert-2004}).
  In this note we show that all geometrically ergodic chains possess
  dominated Coupling from The Past algorithms (not necessarily
  practical!) which are rather closely connected to Foster-Lyapunov
  criteria.
\end{abstract}

\section{Introduction}
\label{sec:intro}

Throughout this paper \(X\) will denote an aperiodic Harris-recurrent
Markov chain on a measurable state space \(\mathcal{X}\) which is a
Polish space (the Polish condition is required in order to
ensure existence of regular conditional probabilities). Recall that \(X\) is
said to be \emph{geometrically ergodic} if it converges in total
variation and at geometric rate to statistical equilibrium \(\pi\), with
multiplicative constant depending on the starting point:
\begin{equation}
  \label{eq:geometric-equilibrium}
\dist_{\text{TV}}(\Law{X_n},\pi) \quad\leq\quad V(X_0) \gamma^n
\end{equation}
for some function \(V:\mathcal{X}\to[1,\infty)\) and some rate \(\gamma\in(0,1)\).
The chain \(X\) is said to be \emph{uniformly ergodic} if the function
\(V\) can be chosen to be constant.

We also recall the notion of a small set:
\begin{definition}\label{lem:minorization}
  A subset \(C\subseteq\mathcal{X}\) is a \emph{small set (of order
  \(k\))} for the Markov chain \(X\) if there is a \emph{minorization
  condition}: for \(\beta\in(0,1)\), and probability measure \(\nu\),
\begin{equation}
  \label{eq:minorization}
  \Prob{X_{k}\in E\;|\;X_0=x} \quad\geq\quad \beta \Indicator{x\in C}\times\nu(E)
\quad\text{for all measurable }E\subseteq\mathcal{X}\,.
\end{equation}
\end{definition}
Results are often stated in terms of the more general notion of
\emph{petite sets}; however for \(\psi\)-irreducible aperiodic chains the
two notions are equivalent \\
\cite[Theorem 5.5.7]{MeynTweedie-1993}.

\citeN{FossTweedie-1997} use small set theory to show that the
condition of uniform ergodicity for such \(X\) is \emph{equivalent} to
the existence of a Coupling from the Past algorithm in the sense of
\citeN{ProppWilson-1996}. This {\emph{classic} \CFTP{}} algorithm
delivers a perfect sample from the equilibrium distribution of \(X\).
The key to the \citeANP{FossTweedie-1997} argument is to remark that
in case of uniform ergodicity the entire state space is small.
Sub-sampling the process \(X\) if necessary (to reduce the
\hyperref[lem:minorization]{order of the small set} to \(1\)), one can
then devise a classic \CFTP algorithm which is actually of the form
introduced by \citeN{MurdochGreen-1997} as the \emph{multigamma
  coupler}.  \citeN{HobertRobert-2004} develop the
\citeANP{FossTweedie-1997} argument to produce approximations to deal
with \emph{burn-in} (time till approximate equilibrium) in the
geometrically ergodic case.

The \citeANP{FossTweedie-1997} result might be thought to delimit and
constrain the possible range of applicability of \CFTP. However it is
also possible to sample perfectly from the equilibrium of some
strictly geometrically ergodic chains using a generalization: namely
\emph{dominated \CFTP} (\domCFTP) as introduced in \citeN{Kendall-1998a},
\citeN{KendallMoller-2000}, \citeN{CaiKendall-1999a}. In this note we
show that this is generic: geometric ergodicity implies the existence
of a special form of \domCFTP algorithm adapted to the geometric
ergodicity in question.  Recent expositions of quantitative
convergence rate estimation depend heavily on small sets and their
relatives (see for example \citeNP{Rosenthal-2002}), so this
piece of \CFTP theory connects to quantitative convergence theory in a
rather satisfying way.

To describe this special form of \domCFTP, we must first introduce the
notion of a Foster-Lyapunov condition.  Geometric ergodicity for our
\(X\) is equivalent to a \emph{geometric Foster-Lyapunov condition}
involving recurrence on small sets (this can be extracted from
\citeNP[Theorem 16.0.1]{MeynTweedie-1993}):
\begin{equation}\label{eq:foster-lyapunov}
\Expect{\Lambda(X_{n+1}) \;|\; X_{n}=x} \quad\leq\quad
\alpha \Lambda(x) + b \Indicator{X_n\in C}\,,
\end{equation}
for some \(\alpha\in(0,1)\) and \(b>0\), some 
\hyperref[lem:minorization]{small set} \(C\), and a function
\(\Lambda:\mathcal{X}\to[1,\infty)\) which is bounded on \(C\).
Note that \(\alpha+b\geq1\) is required, as is \(\Lambda|_{C^c}\geq\alpha^{-1}\), 
since we impose \(\Lambda\geq1\). 
\label{page:marker2} 

Now \hyperref[eq:foster-lyapunov]{Condition
  (\ref*{eq:foster-lyapunov})} implies that every sub-level set
\(\{x\in\mathcal{X}:\Lambda(x)\leq c\}\) is small (as indeed do weaker
  conditions; \citeNP[Theorem 14.2.3]{MeynTweedie-1993}). 
 \newline
This is a key fact for our argument so we
        sketch a coupling proof.

First note that without loss of generality we
       can employ sub-sampling to ensure that the small set \(C\) in
       \hyperref[eq:foster-lyapunov]{Condition
         (\ref*{eq:foster-lyapunov})} is of
       \hyperref[lem:minorization]{order \(1\)}. Super-martingale
       arguments show that we can choose \(n\) such that
       \(\Prob{X\text{ hits }C\text{ before }n\;|\;X_0=x}\) can be
       bounded away from zero uniformly in \(x\) for \(\Lambda(x)\leq c\). Let
       the hitting probability lower bound be \(\rho_0\). We can use the
       \hyperref[eq:minorization]{Minorization Condition
         (\ref*{eq:minorization})} to realize \(X\) as a split-chain in
       the sense of \citeN{Nummelin-1978}, regenerating with
       probability \(\beta\) whenever \(X\in C\).  Couple chains from
       different starting points according to the time when \(X\) first
       regenerates in \(C\), yielding a family of realizations \(X^x\)
       of the Markov chain, with \(X^x_0=x\), such that with positive
       probability \(\beta\rho_0\) all realizations \(\{X^x : \Lambda(x)\leq c\}\)
       coalesce into a set of at most \(n\) trajectories by time \(n\)
       (divided according to the time of first regeneration).  Now
       apply a renewal-theoretic argument to the subsequent
       regenerations of this finite set of trajectories, which are
       allowed to evolve independently, except that whenever two
       trajectories regenerate at the same time they are forced to
       coalesce.  Straightforward analysis shows that we can choose
       \(m\) such that with positive probability \(\rho_1<\beta\rho_0\) all
       trajectories starting from \(\{x\in\mathcal{X}:\Lambda(x)\leq c\}\) have
       coalesced to just one trajectory by time \(n+m\).  Hence
       \(\{x\in\mathcal{X}:\Lambda(x)\leq c\}\) is a small set of order \(n+m\),
       with minorization probability \(\rho_1\).  
%
It is convenient to isolate the notion of a \emph{scale function}
such as \(\Lambda\) in \hyperref[eq:foster-lyapunov]{Equation
  (\ref*{eq:foster-lyapunov})}.
\begin{definition}\label{def:scale-function}
  A \emph{(Foster-Lyapunov) scale function} for a Markov chain state
  space \(\mathcal{X}\) is a measurable function
\[
\Lambda:\mathcal{X}\to[1,\infty)
\]
such that sub-level sets \(\{x\in\mathcal{X}:\Lambda(x)\leq \lambda\}\) are small for all
\(\lambda\geq1\).
\end{definition}

Now we can define the special form of \domCFTP which we require, which
is adapted to a specified Foster-Lyapunov scale function.
\begin{definition}\label{defn:dominating-scale}
  Suppose that \(\Lambda\) is a scale function for an Harris-recurrent Markov chain
  \(X\).  We say the stationary ergodic random process \(Y\) on
  \([1,\infty)\) is a \emph{dominating process for \(X\) based on the scale
    function \(\Lambda\)} (with \emph{threshold} \(h\) and \emph{coalescence
    probability} \(\varepsilon\)) if it is coupled co-adaptively to realizations
of \(X^{x,-t}\) (the Markov chain \(X\) begun at \(x\) at time \(-t\))
as follows:
\begin{itemize}
\item[(a)] for all \(x\in\mathcal{X}\), \(n>0\), and \(-t\leq0\),
almost surely
\begin{equation}\label{eq:domination-requirement}
\Lambda(X^{x,-t}_{-t+n})\quad\leq\quad Y_{-t+n} \qquad\Rightarrow\qquad
\Lambda(X^{x,-t}_{-t+n+1})\quad\leq\quad
Y_{-t+n+1}\,;
\end{equation}
\item[(b)] moreover if \(Y_n\leq h\) then the probability  of
  \emph{coalescence} is at least \(\varepsilon\), where coalescence means that the set 
\[
\left\{ X^{x,-t}_{n+1}\;:\; 
\text{ such that }-t\leq n \text{ and } \Lambda(X^{x,-t}_n)\leq Y_n \right\} 
\]
is a singleton set;
\item[(c)] and finally, \(\Prob{Y_n\leq h}\) must be positive.
\end{itemize}
\end{definition}

Suppose \(Y\) is a dominating process for \(X\) based on the scale
\(\Lambda\).  The following \domCFTP algorithm then yields a draw from
the equilibrium distribution of \(X\).
\begin{algorithm}\label{ag:dom-cftp}
\begin{itemize}
\item[]
\item[] Simulate \(Y\) backwards in equilibrium till the most recent
    \(T<0\) for which \(Y_T\leq h\);
\item[] while coalescence does not occur
    at time \(T\):
  \begin{itemize}
  \item[] extend \(Y\) backwards till the most recent
      \(S<T\) for which \(Y_S\leq h\);
  \item[] set \(T\gets S\);
  \end{itemize}
\item[] simulate the coupled \(X\) forwards from time \(T+1\),
    starting with the unique state produced by the
    coalescence event at time \(T\);
\item[] return \(X_0\) as a perfect draw from equilibrium.
\end{itemize}
\end{algorithm}
Practical implementation considerations are: (1) can one draw from the
equilibrium of \(Y\)? (2) can one simulate \(Y\) backwards in
equilibrium? (3) can one couple the dominated target processes
\(X^{x,-t}\) with \(Y\) so as to ensure the possibility of
regeneration? (4) can one determine when this regeneration has
occurred? and, of course, (5) will the algorithm not run too slowly?

The simplest kind of ordinary small-set \CFTP, as in
\citeN{MurdochGreen-1997}, is recovered from this Algorithm by taking
\(Y\equiv h\), and requiring the whole state-space to be small.  In actual
constructions, care must be taken to ensure that \(Y\) dominates a
coupled collection of \(X\) for which coalescence is possible as
specified in \hyperref[defn:dominating-scale]{Definition
  \ref*{defn:dominating-scale}(b)} (see the treatment of \CFTP for Harris
chains in \citeNP{CorcoranTweedie-2000}).

The proof that this algorithm returns a perfect draw from the
equilibrium distribution of \(X\) is an easy variation on the usual
\domCFTP argument, found at varying levels of generality in
\citeNP{Kendall-1998a,KendallMoller-2000,CaiKendall-1999a}. The key is
to observe that \hyperref[ag:dom-cftp]{Algorithm \ref*{ag:dom-cftp}}
reconstructs a coalesced trajectory which may be viewed as produced by
the Markov chain begun at time \(-\infty\) at some specified state \(x\)
such that \(\Lambda(x)\leq h\): the proof is then an exercise in making this
heuristic precise.

The \citeN{FossTweedie-1997} argument, and the fact that the
geometric Foster-Lyapunov 
\hyperref[eq:foster-lyapunov]{condition
  (\ref*{eq:foster-lyapunov})} would certainly produce a dominating
process if the expectation inequality was replaced by a stochastic
domination, suggests our main result, which will be proved in
\hyperref[sec:implication]{Section \ref*{sec:implication}}:
\begin{theorem}\label{thm:geometric-domCFTP}
  If \(X\) is a geometrically ergodic Markov chain, and \(\Lambda\) is a
  scale function for \(X\) which is derived from some
  geometric Foster-Lyapunov condition, then there exists a \domCFTP
  algorithm for \(X\) (possible subject to sub-sampling) using a
  dominating process based on the scale \(\Lambda\), as in
  \hyperref[ag:dom-cftp]{Algorithm \ref*{ag:dom-cftp}}.
\end{theorem}

As in the case of the \citeN{FossTweedie-1997} result, this algorithm
need not be at all practical!

\section{Geometric ergodicity implies \domCFTP}
\label{sec:implication}
We begin with a lemma concerning the effect of sub-sampling on the 
geometric Foster-Lyapunov
\hyperref[eq:foster-lyapunov]{
  condition}.
\begin{lemma}\label{lem:sub-sampling}
  Suppose \(X\) satisfies a \hyperref[eq:foster-lyapunov]{geometric
    Foster-Lyapunov condition}: for some \(\alpha<1\), some scale function
    \(\Lambda\), and small set \(C=\{x\in\mathcal{X}:\Lambda(x)\leq c\}\).
\begin{equation}\label{eq:fl-working}
\Expect{\Lambda(X_{n+1}) \;|\; X_{n}=x} \quad\leq\quad \alpha \Lambda(x) + b
\Indicator{\Lambda(X_n)\leq c}\,.
\end{equation}
Under \(k\)-sub-sampling we obtain a similar condition but with
different constants:
\begin{equation}
\Expect{\Lambda(X_{n+k}) \;|\; X_{n}=x} \quad\leq\quad \alpha^{k-1} \Lambda(x) + b^{\prime}
\Indicator{\Lambda(X_n)\leq c^{\prime}}\,,
\label{eq:sub-sampling-1}
\end{equation}
and also, if \(k\geq2\),
\begin{equation}
\Expect{\Lambda(X_{n+k}) \;|\; X_{n}=x} \quad\leq\quad \alpha \Lambda(x) + b^{\prime\prime}
\Indicator{\Lambda(X_n)\leq c^{\prime\prime}}\,.
\label{eq:sub-sampling-2}
\end{equation}
Moreover \(b^\prime=b/(1-\alpha)\), \(c\prime=b/(\alpha^{k-1} (1-\alpha)^2)\) may be chosen not
to depend on \(c\), and \(b^{\prime\prime}=b/(1-\alpha)\), \(c^{\prime\prime}=b/(\alpha(1-\alpha)^2)\)
may be chosen to depend neither on \(c\) nor on \(k\geq2\).
\end{lemma}
We are able to choose \(b^\prime\), \(c\prime\), \(b^{\prime\prime}\), \(c^{\prime\prime}\) not to
depend on \(c\) because we have allowed generous sub-sampling
(\emph{i.e.}: \(k\)-sub-sampling to change \(\alpha\) to \(\alpha^{k-1}\)).

\begin{proof}
  Iterating \hyperref[eq:fl-working]{Equation (\ref*{eq:fl-working})},
  \begin{align*}
\Expect{\Lambda(X_{n+k}) \;|\; X_{n}=x} &\quad\leq\quad    
\alpha^k \Lambda(x) + \sum_{j=1}^k \alpha^{j-1} b\Expect{\Indicator{\Lambda(X_{n+k-j})\leq c}\;|\;X_{n}=x}\\
&\quad\leq\quad    
\alpha^k \Lambda(x) + \frac{b}{1-\alpha} \\
&\quad=\quad
\alpha^{k-1} \Lambda(x) - \alpha^{k-1} (1-\alpha) \Lambda(x)+ \frac{b}{1-\alpha} \\
&\quad\leq\quad
\begin{cases}
  \alpha^{k-1} \Lambda(x) & \text{ if } \Lambda(x)> \frac{b}{\alpha^{k-1} (1-\alpha)^2} \,,\\
  \alpha^{k-1} \Lambda(x) + {b}/{(1-\alpha)} & \text{ otherwise.} 
\end{cases}
  \end{align*}
  Hence we may choose \(b^{\prime}=b/(1-\alpha)\), \(c^{\prime}=b/(\alpha^{k-1}
  (1-\alpha)^2)\). Alternatively
\begin{align*}
  \Expect{\Lambda(X_{n+k}) \;|\; X_{n}=x} &\quad\leq\quad  
\alpha \Lambda(x) - \alpha (1-\alpha^{k-1}) \Lambda(x)+ \frac{b}{1-\alpha} \\
&\quad\leq\quad
\begin{cases}
  \alpha \Lambda(x) & \text{ if } \Lambda(x)> \frac{b}{\alpha (1-\alpha)(1-\alpha^{k-1})} \,,\\
  \alpha \Lambda(x) + {b}/{(1-\alpha)} & \text{ otherwise.} 
\end{cases}
\end{align*}
Hence we may choose \(b^{\prime\prime}=b/(1-\alpha)\), \(c^{\prime\prime}=b/(\alpha (1-\alpha)^2)\) if
\(k\geq2\).
\end{proof}

\begin{proof}[of Theorem \ref*{thm:geometric-domCFTP}]
We first construct the dominating process.

Consider Markov's inequality applied to the
geometric Foster-Lyapunov 
\hyperref[eq:foster-lyapunov]{inequality
  (\ref*{eq:foster-lyapunov})}.
Any dominating process \(Y\) must satisfy the
\hyperref[eq:domination-requirement]{stochastic domination
  (\ref*{eq:domination-requirement})} described in
\hyperref[defn:dominating-scale]{Definition
  \ref*{defn:dominating-scale}}. Consequently, in default of further
distributional information about \(\Prob{\Lambda(X_{n+1}) | X_{n}=x}\),
if \(Y\) is to
be a dominating process based on the scale \(\Lambda\) then we need \(Y\) to
be stationary ergodic but also to satisfy
\begin{equation}\label{eq:crude-domination}
   \Prob{Y_{n+1}\geq \alpha z y\;|\; Y_n=z} \quad\geq\quad 
 \sup_{x: \Lambda(x)\leq z}\frac{\Expect{\Lambda(X_{n+1})\;|\;X_n=x}}{\alpha z y}\,.
\end{equation}

Now if \(C\subseteq\{x\in\mathcal{X}:\Lambda(x)\leq c\}\) then
 \begin{align*}
 \sup_{x: \Lambda(x)\leq z}\frac{\Expect{\Lambda(X_{n+1})\;|\;X_n=x}}{\alpha z y}
 &\quad\leq\quad
 \sup_{x: \Lambda(x)\leq z} 
\frac{\alpha\Lambda(x) + b \Indicator{x:\Lambda(x)\leq c}}{\alpha z y}\\
 \quad\leq\quad \sup_{x: \Lambda(x)\leq z}\frac{\alpha \Lambda(x)}{\alpha z y} & \quad=\quad \frac{1}{y}
 \qquad\text{ so long as }
 z\geq c + \frac{b}{\alpha}\,.
 \end{align*}

Consequently \(Y\) is a possible candidate for a dominating process based
on the scale \(\Lambda\) if
\begin{equation}
  \label{eq:dominating-in-scale}
  \Prob{Y_{n+1}\geq \alpha z y \:|\; Y_n=z} \quad=\quad
  \begin{cases}
    1/y & \text{ if } z\geq c + \frac{b}{\alpha} \,,\\
    1   & \text{ otherwise.}
  \end{cases}
\end{equation}
If we define \(U\) by \(Y=(c+b/ \alpha)\exp(U)\) (so \(U\) is a
\emph{log-dominating process}) then \(U\) is the system workload of a
\(D/M/1\) queue, sampled at arrivals, with arrivals every \(\log(1/
\alpha)\) units of time, and service times being independent and of unit
Exponential distribution. The process \(U\) is a random walk with
reflection (of Skorokhod type) at \(0\): as its jump distribution is
\(\text{Exponential}(1)-\log(1/ \alpha)\) we may deduce it is
positive-recurrent if and only if \(\alpha<e^{-1}\).

\label{page:marker1}
In case \(e^{-1}<\alpha<1\), \(U\) and \(Y=(c+b/ \alpha)\exp(U)\) fail to be
positive-recurrent.  However the same construction will work if we use
\hyperref[eq:sub-sampling-1]{Equation (\ref*{eq:sub-sampling-1})} of
\hyperref[lem:sub-sampling]{Lemma \ref*{lem:sub-sampling}} to justify
sub-sampling \(X\) with a sampling period \(k\) large enough to ensure
a \hyperref[eq:foster-lyapunov]{geometric Foster-Lyapunov condition
  (\ref*{eq:foster-lyapunov})} using \(\Lambda\) as scale but with \(\alpha\)
replaced by \(\alpha^{k-1}<e^{-1}\), and amending \(b\) to \(b^\prime\), \(c\)
to \(c^\prime\) as in \hyperref[eq:sub-sampling-1]{Inequality
  (\ref*{eq:sub-sampling-1})}.

Thus without loss of generality 
we may assume \(\alpha<e^{-1}\), 
and so this \(Y\) can be run in statistical equilibrium, and thus
qualifies as least partly as a dominating process for the purposes of
\hyperref[thm:geometric-domCFTP]{Theorem
  \ref*{thm:geometric-domCFTP}}. In the sequel we assume moreover that
further sub-sampling has been carried out based on
\hyperref[eq:sub-sampling-2]{Equation (\ref*{eq:sub-sampling-2})}, to
ensure that the following small set is of order \(1\):
\begin{equation}
  \label{eq:small-target}
  \left\{x\in\mathcal{X}\;:\; \Lambda(x) \leq h \right\}
\qquad\text{ for }\qquad
h=\max\left\{c+\frac{b}{\alpha},
\frac{b}{\alpha(1-\alpha)}\left(1+\frac{1}{1-\alpha}\right)
\right\}\,.
\end{equation}
Here the level \(h\geq c+b/ \alpha\) is fixed so as to ensure \(h=c^{\prime\prime}+b^{\prime\prime}/(1-\alpha)\) with
\(b^{\prime\prime}\), \(c^{\prime\prime}\) given as in
\hyperref[eq:sub-sampling-2]{Equation (\ref*{eq:sub-sampling-2})};
thus \(h\) supplies a stable threshold for geometric Foster-Lyapunov
conditions, even allowing for further sub-sampling if required. Note
in particular that \(Y=(c+b/
\alpha)\exp(U)\) is able to sink below \(h\), since \(h\geq c+b/\alpha\) and the
system workload \(U\) can reach zero.
\label{page:marker3}

To fulfil the requirements on a dominating process given in
\hyperref[defn:dominating-scale]{Definition
  \ref*{defn:dominating-scale}}, we need to construct a coupling
between \(Y\) and the target process \(X\) expressed in
terms of a random flow of independent maps
\(F_{-t+n+1}:\mathcal{X}\to\mathcal{X}\):
\[
X^{x,-t}_{-t+n+1}\quad=\quad F_{-t+n+1}(X^{x,-t}_{-t+n})
\]
satisfying the distributional requirement that \(X^{x,-t}\) should
evolve as the Markov chain \(X\), 
the \hyperref[eq:domination-requirement]{domination requirement
expressed by the implication (\ref*{eq:domination-requirement})}, 
and also the regeneration requirement
that with probability \(\varepsilon\) the set
\[
\left\{ F_n(u) \;:\; \text{ such that } \Lambda(u)\leq h \right\} 
\]
should be a singleton set. The well-known link between stochastic
domination and coupling can be applied together with the arguments
preceding \hyperref[eq:dominating-in-scale]{Equation
  (\ref*{eq:dominating-in-scale})} to show that we can couple the
various \(X^{x,-t}\) with \(Y\) co-adaptively in this manner so that
the implication (\ref*{eq:domination-requirement}) holds: note that
here and here alone we use the Polish space nature of \(\mathcal{X}\),
which allows us to complete the couplings by constructing regular
conditional probability distributions for the various \(X^{x,-t}\)
conditioned on the \(\Lambda(X^{x,-t})\). Thus all that is required is to
show that this stochastic domination coupling can be modified to allow
for regeneration.

The small set condition for \(\{x\in\mathcal{X}:\Lambda(x)\leq h\}\) means there is
a probability measure \(\nu\) and a scalar \(\beta\in(0,1)\) such that for all
Borel sets \(B\subseteq[1,\infty)\), whenever \(\Lambda(x)\leq h\),
\begin{equation}
  \label{eq:crucial-small-set}
  \Prob{\Lambda(X_{n+1})\in B \;|\; X_n=x} \quad\geq\quad \beta\nu(B) \,.
\end{equation}
Moreover the stochastic domination which has been arranged in the
 course of defining \(Y\) means that for all real \(u\),
whenever \(\Lambda(x) \leq y\),
 \begin{equation}
   \label{eq:crucial-stochastic-domination}
  \Prob{\Lambda(X_{n+1})>u \;|\; X_n=x} \quad\leq\quad 
\Prob{Y>u\;|\; Y=y}\,.
 \end{equation}
We can couple in order to arrange for regeneration if we can identify
a probability measure \(\widetilde\nu\), defined solely in terms of
\(\nu\) and the dominating jump distribution \(\Prob{Y\geq u \;|\;
  Y=y}\), such that for all real \(u\)
\begin{align*}
\Prob{\Lambda(X_{n+1})>u \;|\; X_n=x} - \beta \nu((u,\infty))
\quad  &\leq\quad\Prob{Y>u \;|\;  Y=y} - \beta\widetilde\nu((u,\infty))\\
\nu((u,\infty)) \quad &\leq\quad \widetilde\nu((u,\infty))
\end{align*}
and moreover
\[
  \Prob{Y_{n+1}\in B \;|\; Y_n=y} \quad\geq\quad \beta\widetilde\nu(B) \,.
\]
For then at each step we may determine whether or not regeneration has
occurred (with probability \(\beta\)); under regeneration we use
stochastic domination to couple \(\nu\)
to \(\widetilde\nu\); otherwise we use
stochastic domination to couple the residuals.

We state and prove this as an interior lemma, as it may be of wider interest.
\begin{lemma}\label{lem:mixture-domination-coupling}
Suppose \(U\), \(V\) are two random variables defined on \([1,\infty)\) such
that
\begin{itemize}
\item[(a)] The distribution \(\Law{U}\) is stochastically dominated by the distribution 
  \(\Law{V}\):
\begin{equation}
\Prob{U> u} \quad\leq\quad \Prob{V>u}\qquad\text{ for all real }U\,;
\label{eq:lemma-domination}
\end{equation}
\item[(b)] \(U\) satisfies a minorization condition: for some
  \(\beta\in(0,1)\) and probability measure \(\nu\):
  \(B\subseteq[1,\infty)\),
\begin{equation}
\Prob{U\in B} \quad\geq\quad \beta\nu(B) \qquad \text{ for all Borel sets }B\subseteq[1,\infty)\,.
\label{eq:lemma-minorization}
\end{equation}
\end{itemize}
Then there is a probability measure \(\mu\) stochastically
dominating \(\nu\) and such that \(\beta\mu\) is minorized by
\(\Law{V}\). Moreover \(\mu\) depends only on \(\beta\nu\)
and \(\Law{V}\).
\end{lemma}
\begin{proof}[of Lemma \ref{lem:mixture-domination-coupling}]
  Subtract the measure \(\beta\nu((u,\infty))\) from both sides of
  \hyperref[eq:lemma-domination]{Inequality (\ref*{eq:lemma-domination})}
  representing the stochastic domination \(\Law{U}\preceq\Law{V}\). By the
\hyperref[eq:lemma-minorization]{minorization condition
  (\ref*{eq:lemma-minorization})} the resulting left-hand-side is nonnegtive. Thus for all real
  \(u\)
\[
0 \quad\leq\quad \Prob{U>u} - \beta\nu((u,\infty)) \quad\leq\quad \Prob{V>u} - \beta\nu((u,\infty))
\]
 Now
  \(\Law{U}-\beta\nu\) is a nonnegative measure (because of the
  \hyperref[eq:lemma-minorization]{minorization condition
  (\ref*{eq:lemma-minorization})}). Consequently \(\Prob{U>u} - \beta\nu((u,\infty))\)  must be
  non-increasing in \(u\) and so we may reduce the
  right-hand side by minimizing over \(w\leq u\):
\begin{align*}
  \Prob{U>u} - \beta\nu((u,\infty)) &\quad\leq\quad \inf_{w\leq u}\left\{ \Prob{V>w} -
    \beta\nu((w,\infty)) \right\}\\
&\quad=\quad \Prob{V>u} -\beta\mu((u,\infty))
\end{align*}
where \(\mu\) is the potentially \emph{signed} measure defined by
\[ 
\beta \mu([1,u]) \quad=\quad
\Prob{V\leq u} - \sup_{w\leq u}\left\{ \Prob{V\leq w} - \beta\nu([1,w)) \right\}\,.
\]
In fact \(\mu\) is a probability measure on \([1,\infty)\). Both
\(\mu(\{1\})=\nu(\{1\})\) and
\(\mu([1,\infty))=1\) follow from considering \(u=1\), \(u\to\infty\). Now we show
\(\mu\) is nonnegative:
\begin{align*}
&  \beta\mu((u,u+u^\prime]) - \Prob{u<V\leq u+u^\prime} 
 \\ 
& \quad=\quad
 - \sup_{w\leq u+u^\prime}\left\{ \Prob{V\leq w} - \beta\nu([1,w)) \right\}  
 + \sup_{w\leq u}\left\{ \Prob{V\leq w} - \beta\nu([1,w)) \right\}\,.
\end{align*}
If the first supremum were to be attained at \(w\leq u\) then the two suprema
would cancel. If the first supremum were to be attained at \(w^\prime\in[u,u+u^\prime]\)
then
\begin{align*}
&  \beta\mu((u,u+u^\prime]) - \Prob{u<V\leq u+u^\prime} 
 \\
& \quad=\quad
 - \Prob{V\leq w^\prime} + \beta\nu([1,w^\prime))
 + \sup_{w\leq u}\left\{ \Prob{V\leq w} - \beta\nu([1,w)) \right\}\\
& \quad\geq\quad
 - \Prob{V\leq w^\prime} + \beta\nu([1,w^\prime))
 + \Prob{V\leq u} - \beta\nu([1,u)
\end{align*}
and hence
\[
\beta\mu((u,u+u^\prime])
\quad\geq\quad \Prob{w^\prime<V\leq u+u^\prime}
 + \beta\nu([u,w^\prime))
\quad\geq\quad0\,.
\]
So we can deduce \(\beta\mu\) is in fact a nonnegative measure.
\label{page:marker4}

On the other hand
\begin{align*}
&  \beta\mu((u,u+u^\prime]) - \Prob{u<V\leq u+u^\prime} 
 \\
& \quad=\quad
 - \sup_{w\leq u+u^\prime}\left\{ \Prob{V\leq w} - \beta\nu([1,w)) \right\} 
 + \sup_{w\leq u}\left\{ \Prob{V\leq w} - \beta\nu([1,w)) \right\} \\
&\quad\leq\quad0\,,
\end{align*}
hence
\begin{equation}
  \label{eq:sandwich}
   0\quad\leq\quad \beta\mu((u,u+u^\prime])\quad\leq\quad \Prob{u<V\leq u+u^\prime}\,,
\end{equation}
so \(\beta \mu\) is absolutely continuous with respect to
\(\Law{V}\) and indeed we can deduce 
\begin{equation}
  \label{eq:representation}
  \beta\d\mu(u) \quad=\quad \Indicator{\Prob{V>\cdot} -
    \beta\nu((\cdot,\infty))\text{ hits current minimum at }u } \d\Prob{V\leq u}\,.
\end{equation}
The minorization of \(\beta\mu\) by \(\Law{V}\) follows from this argument:
dependence
only on \(\beta\nu\)
and \(\Law{V}\) follows by construction; finally, stochastic
domination of \(\beta \nu\) follows from
\begin{align*}
  \beta\mu((u,\infty))& \quad=\quad \Prob{V>u} - \inf_{w \leq u}\left\{ \Prob{V>w} -
      \beta\nu((w,\infty))\right\} \\
& \quad=\quad\sup_{w\leq u}\left\{\beta\nu((w,\infty)) - \Prob{w<V\leq u} \right\} \\
&\quad\geq\quad \beta\nu((u,\infty))\,.
\end{align*}\end{proof}

This concludes the proof of Theorem \ref*{thm:geometric-domCFTP}:
use \hyperref[lem:mixture-domination-coupling]{Lemma
  \ref*{lem:mixture-domination-coupling}} to couple
\(\Law{X_{n+1}\;|\;X_n=x}\) to \(\Law{Y_{n+1}\;|\;Y_n=y}\) whenever
\(\Lambda(x)\leq y\) in a way which implements stochastic domination and
ensures all the \(X_{n+1}\) regenerate simultaneously whenever \(Y\leq h\).
\end{proof}

Note that the algorithm requires us to be able to draw from the
equilibrium distribution of \(Y\) and to simulate its time-reversed
equilibrium dual. Up to an additive constant \(\log(Y)\) is the workload
of a \(D/M/1\) queue. This queue is amenable to exact calculations, so
these simulation tasks are easy to implement 
(specializing the theory of the \(G/M/1\) queue as discussed
in \citeNP[ch.~11]{GrimmettStirzaker-1992}). However in general we do \emph{not} expect this
``universal dominating process'' to lead to practical \domCFTP
algorithms! The difficulty in application will arise in determining whether or
not regeneration has occurred as in \hyperref[ag:dom-cftp]{Algorithm
  \ref*{ag:dom-cftp}}. This will be difficult especially if
sub-sampling has been applied, since then one will need detailed
knowledge of convolutions of the probability kernel for \(X\)
(potentially a harder problem than sampling from equilibrium!).

Of course, in practice one uses different dominating processes
better adapted to the problem at hand. For example an \(M/D/1\) queue
serves as a good log-dominating process for perpetuity-type problems
and gives very rapid \domCFTP algorithms indeed, especially when
combined with other perfect simulation ideas such as multishift \CFTP
\cite{Wilson-1999}, read-once \CFTP \cite{Wilson-2000a}, or one-shot
coupling \cite{RobertsRosenthal-2002}.

Finally note that, in cases when \(\alpha\in[e^{-1},1)\) or when the small
set \(\{x\in\mathcal{X}:\Lambda(x)\leq h\}\) is of order greater than \(1\), we are
forced to work with coupling constructions that are effectively
\emph{non-co-adapted} (sub-sampling means that target transitions
\(X_{mk}\) to \(X_{mk+1}\) depend on sequences \(Y_{mk}\), \(Y_{mk+1}\), \ldots,
\(Y_{mk+k}\)). The potential improvements gained by working with
non-adapted couplings are already known not only to theory (the
non-co-adapted filling couplings of
\citeNP{Griffeath-1974,Goldstein-1978}; and the efficiency
considerations of \citeNP{BurdzyKendall-1997a}) but also to
practitioners (\citeNP{Huber-2004}: non-Markovian techniques in \CFTP;\\
\citeNP{HayesVigoda-2003}: non-Markovian conventional MCMC for random
sampling of colorings).

\section{Counter-example}
\label{sec:counter-example}
We complete this note by describing a counter-example:
a Markov chain \(X\) which satisfies a Foster-Lyapunov
condition involving a scale function \(\Lambda\), but such that there can be no
recurrent dominating process \(Y\) based on \(\Lambda\).  We begin by
choosing a sequence of disjoint measurable sets \(S_1\), \(S_2\), \ldots,
subsets of \([1,\infty)\) such that each set places positive measure in
every non-empty open set:
\begin{lemma}\label{lem:partition}
  One can construct a measurable partition \(S_1\), \(S_2\), \ldots of \([1,\infty)\),
\[
S_1 \sqcup S_2 \sqcup S_3 \sqcup \ldots \quad=\quad [1,\infty)\,,
\]
with the property
\(\Leb(S_i \cap (u,v)) > 0\)
for all \(0<u<v<\infty\), all \(i\in\{1, 2, \ldots\}\).
\end{lemma}
\begin{proof}
  Enumerate the rational numbers in \([0,1)\) 
by \(0=\tilde q_0\), \(\tilde q_1\), \(\tilde q_2\), \ldots . Choose \(\alpha<1/2\), and define
\[
A_0 \quad=\quad \bigcup_{k=1}^\infty \bigcup_{n=0}^\infty \left[\tilde q_n+k,\tilde q_n+k+\alpha
  2^{-n}\right]\,.
\]
Then for each \(k\geq1\)
\[
\alpha \quad\leq\quad \Leb\left(A_0\cap[k,k+1)\right) \quad\leq\quad 2\alpha\,.
\]
Continue by defining a sequence of nested subsets \(A_{r}\subset A_{r-1}\) by
\begin{equation}
  \label{eq:subsets}
A_r \quad=\quad 
\bigcup_{k=1}^\infty \bigcup_{n=0}^\infty \left[\frac{\tilde q_n+k}{2^r},\frac{\tilde q_n+k}{2^r}+\frac{\alpha}{4^r} 2^{-n}\right]\,,
\end{equation}
satisfying
\begin{equation}
  \label{eq:subsets-bounds}
\frac{\alpha}{4^r} \quad\leq\quad \Leb\left(A_r\cap\Big[\frac{k}{2^r},\frac{k+1}{2^r}\Big)\right) 
\quad\leq\quad \frac{2\alpha}{4^r}\,.
\end{equation}

Thus the measurable shell \(B_r=A_r\setminus A_{r+1}\) places mass of at least
\(\frac{\alpha}{2\times4^{r}}\) in each interval
\([\frac{k}{2^r},\frac{k+1}{2^r})\)\,.

It follows that if \(S\) is defined by
\[
S \quad=\quad \bigcup_{s=1}^\infty \left(A_{r_s}\setminus A_{r_{s}+1}\right)
\]
then \(\Leb(S\cap U)>0\) for every open set \(U\subset[1,\infty)\). The desired
disjoint sequence \(S_1\), \(S_2\), \ldots is obtained by considering a
countably infinite family of disjoint increasing subsequences of the
natural numbers.
\end{proof}

\begin{lemma}
  There is a Markov chain \(X\) satisfying a Foster-Lyapunov condition
  with scale function \(\Lambda\), such that any
  dominating process \(Y\) based on \(\Lambda\) will fail to be positive-recurrent.
\end{lemma}

\begin{proof}
  The Markov chain \(X\) will have state space \([1,\infty)\), with scale
  function \(\Lambda(x)\equiv x\).  We begin by fixing \(\alpha\in(e^{-1},1)\), and set
  \(C=[1,\alpha^{-1}]\). The set \(C\) will be the small set for the
  Foster-Lyapunov condition.  Choose a measurable partition \(S_1 \sqcup
  S_2 \sqcup S_3 \sqcup \ldots = [1,\infty)\) as in \hyperref[lem:partition]{Lemma
    \ref*{lem:partition}}. Enumerate the rational numbers in \([1,\infty)\)
  by \(q_1\), \(q_2\), \ldots.

We define the transition kernel \(p(x,\cdot)\) of \(X\) on
  \([1,\infty)\) as follows:
\begin{itemize}
\item[] For \(x\in[1,\alpha^{-1}]\), set
\[
p(x,\d y) \quad=\quad \exp(-(y-1)) \d y \quad\text{for }y\geq1\,,
\]
so that if \(X_n\in C\) then \(X_{n+1}-1\) has a unit rate
Exponential distribution. Then:
\begin{itemize}
\item[] \(C\) is a small set for \(X\) of order \(1\) (in fact it
  will be a regenerative atom!);
\item[] if \(X_n\in C\) then \(\Expect{X_{n+1}}=2\);
\item[] if \(X\) has positive chance of visiting state \(1\) then the
  whole state space \([1,\infty)\) will be maximally \(\Leb\)-irreducible.
\end{itemize}
\item[] For \(x>\alpha^{-1}\) and \(x\in S_i\), set
\[
p(x, \d y) \quad=\quad
\left(1-\frac{\alpha}{q_i}\right)\delta_{0}(\d y) + \frac{\alpha}{q_i}\delta_{q_i x}(\d y)\,.
\]
Note that, because we are using the identity scale \(\Lambda(x)\equiv x\),
\begin{itemize}
\item[] if \(x\not\in C\) then
  \(\Expect{\Lambda(X_{n+1})\;|\;X_n=x}=\Expect{X_{n+1}\;|\;X_n=x}=\alpha x\);
\item[] if \(x\not\in C\) then \(\Prob{X_{n+1}=1\;|\;X_n=x}>0\).
\end{itemize}
\end{itemize}
Thus \(X\) satisfies a geometric Foster-Lyapunov condition based on
scale \(\Lambda\) and small set \(C\), and so is geometrically ergodic.

Suppose \(Y\) is a dominating process for \(X\) based on the identity
scale \(\Lambda\). This
means it must be possible to couple \(Y\) and \(X\) such that, if
\(\Lambda(X_n)=X_n\leq Y_n\) then \(\Lambda(X_{n+1})=X_{n+1}\leq Y_{n+1}\). This can be achieved if
and only if 
\[
\Prob{X_{n+1}\geq z \;|\; X_n=u} \quad\leq\quad \Prob{Y_{n+1}\geq z\;|\; Y_n=x}\,\
\]
for all \(z\geq1\), and Lebesgue-almost all \(u<x\).
Therefore we require of such 
\(Y\) that
\begin{align*}
&  \Prob{Y_{n+1}\geq \alpha x y\;|\; Y_n=x} \quad\geq\quad 
\text{ess}\sup_{u<x} \left\{ \Prob{X_{n+1}\geq \alpha x y \;|\; X_n=u} \right\}\\
&\quad=\quad
\sup_i \text{ess}\sup\left\{ \frac{\alpha}{q_i}\;:\; \alpha^{-1}<u<x, u\in S_i, q_i
  u > \alpha x y\right\} \\
&\quad=\quad
\sup_i \left\{\frac{\alpha}{q_i}\;:\; q_i > \alpha y \right\} 
\quad=\quad \frac{1}{y}\,,
\end{align*}
using Markov's inequality, then the construction of the kernel of
\(X\), then the measure-density of the \(S_i\).

So such a Markov chain \(Y\) must also (at least when above level
\(\alpha^{-1}\)) dominate \(\exp(Z)\), where \(Z\) is a random walk with
jump distribution \(\text{Exponential}(1)+\log(\alpha)\). Hence it will
fail to be positive-recurrent on the small set \(C\) when \(\alpha\geq e^{-1}\).
\end{proof}

There may exist some subtle re-ordering to provide \domCFTP for such a
chain on a different scale; however the above lemma shows that
\domCFTP must fail for dominating processes for \(X\) based on the
scale \(\Lambda\).

\section{Conclusion}
\label{sec:conclusion}

We have shown that geometric ergodicity (more strictly, a geometric
Foster-\-Lyapunov condition) implies the existence of a special kind of
\domCFTP algorithm. The algorithm is not expected to be practical:
however it connects perfect simulation firmly with more theoretical
convergence results in the spirit of the \citeN{FossTweedie-1997}
equivalence between classic \CFTP and uniform ergodicity. Note also
that the ``universal dominating process'', the sub-critical
\(\exp(D/M/1)\) so derived, is itself geometrically ergodic.

It is natural to ask whether other kinds of ergodicity (for example,
polynomial ergodicity) can also be related to perfect simulation in
this way; this is now being pursued by Stephen Connor as part of his
PhD research at Warwick.

 \bibliographystyle{chicago}

  \thispagestyle{plain}
  \markboth{REFERENCES}{REFERENCES}

 \bibliography{habbrev,ge,wsk}


\end{document}